\documentclass{article}
\usepackage[latin1]{inputenc}
\usepackage{amsfonts}
\usepackage{amsmath, latexsym}
\usepackage{amssymb,verbatim}
\usepackage{theorem}
\usepackage{graphicx, psfrag}
\usepackage{color}
\usepackage{euscript}

\newtheorem{pro}    {Proposition}
\newtheorem{thm}    {Theorem}

\newtheorem{df}     {Definition}

\theorembodyfont{\rmfamily}
\newtheorem{rem}{Remark}

\newcommand{\rv}{\vec{r}}
\newcommand{\uv}{\vec{\theta}}

\newcommand{\wup}{W_{div}^{1/p,p}(\Om, \Ss^1)}

\newcommand{\ba}{\begin{align*}}
\newcommand{\ea}{\end{align*}}

\newcommand{\RR}{\mathbb{R}}
\newcommand{\R}{\mathbb{R}}
\newcommand{\Ss}{\mathbb{S}}

\newcommand{\eps}{\varepsilon}

\newcommand{\h}{{\mc{H}}}

\newcommand{\proof}[1]{\par\medskip\noindent{\bf Proof#1.}}
\newcommand{\qed}{\hfill$\square$}

\newcommand{\be}{\begin{equation}}

\newcommand{\ee}{\end{equation}}

\newcommand{\loc}{_{loc}}

\newcommand{\nd}{\noindent}

\newcommand{\Om}{\Omega}
\newcommand{\dist}{\mathop{\rm dist \,}}

\newcommand{\wun}{W_{div}^{1,1}(\Om, \Ss^1)}

\newcommand{\f}{\varphi}

\newcommand{\supp}{\operatorname{supp}}

\newcommand{\hb}{\mathbf{H}}

\def\I[#1]{\mc{I}_{#1}}

\def\E[#1]{\mc{E}_{#1}}
\def\S[#1]{\mc S_0(#1)}
\def\mc{\mathcal}
\def\ncd{\nabla\cdot}

\def\ed0{\eps\downarrow0}

\title{A regularizing property of the $2D$-eikonal equation}

\author{ {\Large Camillo De Lellis} \footnote{Institut f\"ur Mathematik, Universit\"at Z\"urich, CH-8057 Z\"urich, Switzerland. Email: camillo.delellis@math.uzh.ch} \and
{\Large Radu Ignat}
\footnote{Institut de Math\'ematiques de Toulouse, Universit\'e
Paul Sabatier, 31062
Toulouse, France. Email: Radu.Ignat@math.univ-toulouse.fr}
}

\begin{document}

\maketitle

\begin{abstract}
We prove that any $2$-dimensional solution $\psi\in W\loc^{1+\frac 1 3, 3}$ of the eikonal equation has locally Lipschitz gradient $\nabla \psi$ except at a locally finite number of vortices. 
\end{abstract}

\section{Introduction}
\label{chap_div}

\begin{figure}[ht] %
  \begin{minipage}{0.44\linewidth}
  \centering
   \includegraphics[height=2.5cm]{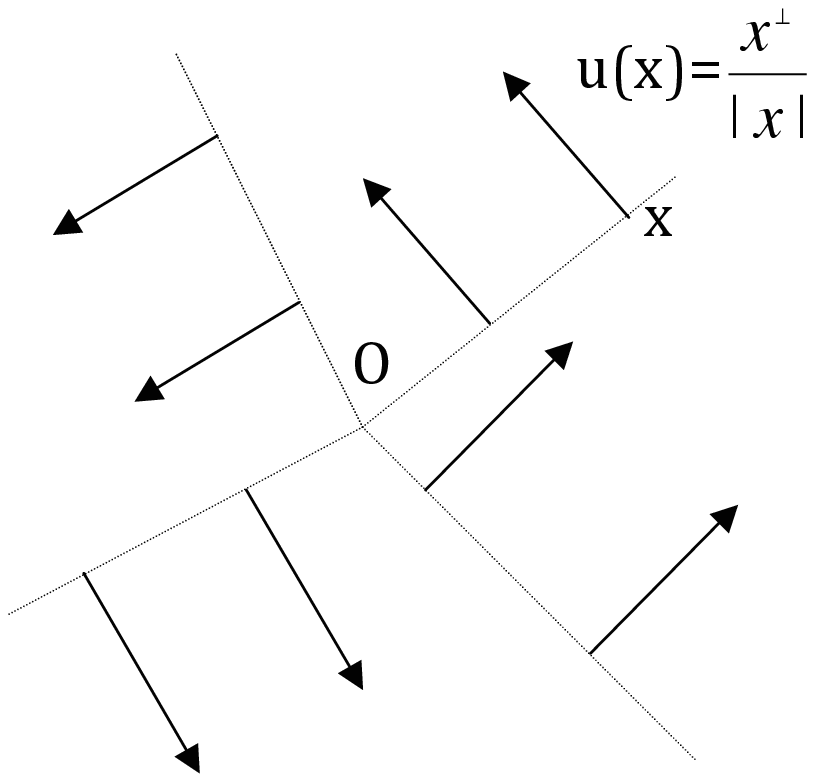} 
     \end{minipage}
  \begin{minipage}{0.52\linewidth}
   \centering
   \includegraphics[height=1.5cm]{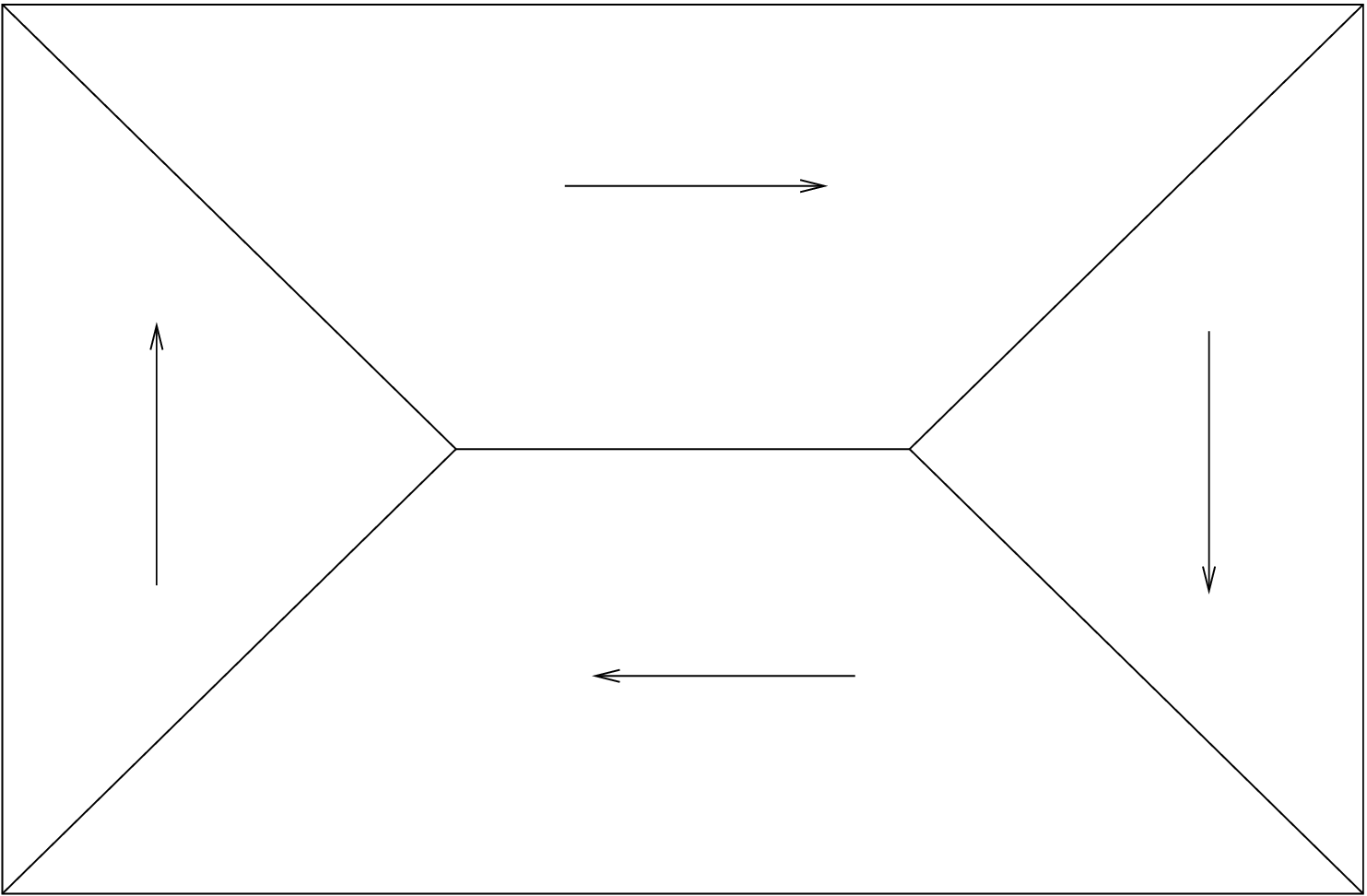} 
  \end{minipage}
   \caption{Vector fields $\nabla^\perp \dist(\cdot, K)$ when $K$ is a point (left) and a rectangle (right).}
  \label{Landau}
\end{figure}

Let $\Omega\subset \RR^2$ be an open set. We will focus on (locally) Lipschitz solutions $\psi:\Omega \to \RR$ of the eikonal equation, namely such that
\be
\label{eiko}
|\nabla \psi|=1 \quad \textrm{a.e. in}\quad  \Omega\, .
\ee
Since all our results will have a local nature, this amounts to investigate curl-free $L^1$ vector fields $w:\Omega\to \RR^2$ of unit length or,
equivalently, $L^1$ vector fields $u=w^\perp=(-w_2, w_1):\Omega\to \RR^2$ that satisfy \be
\label{contraintes}
|u|=1 \, \, \textrm{a.e. in } \Omega \quad \textrm{and} \quad \nabla \cdot u=0\quad \textrm{ in }\, \, {\cal D}'(\Omega).
\ee
Typical examples of stream functions $\psi$ satisfying \eqref{eiko} are distance functions $\psi=\dist(\cdot, K)$ to some closed nonempty set $K\subset \RR^2$ (see Figure \ref{Landau}). 
In general such $\psi$ are not smooth and generate line singularities or vortex-point singularities for the gradient $\nabla \psi$.

We denote by $W_{div}^{s,p} (\Omega, \Ss^1)$  the Sobolev space of order $s>0$ and $p\geq 1$ of divergence-free unit-length vector fields, namely $$W_{div}^{s,p}(\Om, \Ss^1)=\{u\in W\loc^{s,p}(\Om, \RR^2)\, :\, u \textrm{ satisfies } \eqref{contraintes} \, \} ,$$ 
and we show that elements in the critical spaces $u\in W_{div}^{1/p,p} (\Omega, \Ss^1)$ have, for $p\in [1,3]$, only vortex-point singularities, 
i.e. they gain more regularity.

\begin{thm}
\label{teo}
If $u\in \wup$ with $p\in [1, 3]$ then
$u$ is locally Lipschitz continuous inside $\Omega$ except at a locally finite number of singular points. Moreover, every singular point $P$ of $u$ corresponds to a vortex-point singularity of degree $1$ of $u$, i.e., there
exists a sign $\alpha=\pm 1$ such that $$u(x)=\alpha \frac{(x-P)^\perp}{|x-P|} \quad \textrm{for every $x\neq P$ in any convex neighborhood of $P$ in $\Omega$}.$$  
\end{thm}

Following the same strategy we can also show a related regularizing effect for solutions of the Burgers' equation
\begin{equation}\label{e:Burgers}
v_t + \left(\textstyle{\frac{v^2}{2}}\right)_s = 0\, ,
\end{equation}
where $(t,s) = (x_1, x_2)$ will be used for the time-space variables. The link between \eqref{contraintes} and \eqref{e:Burgers} is discussed in the next section. 

\begin{thm}
\label{teo_2}
Let $\Omega = I \times J$ with $I, J\subset \R$ two intervals and $v\in L^4(\Omega)$ be a distributional solution of \eqref{e:Burgers} which belongs to the space $L^3 (I, W^{1/3, 3} (J))$, namely
\begin{equation}\label{e:time-space}
\int_I \int_{J\times J} \frac{|v(t,s) - v (t, \sigma)|^3}{|s-\sigma|^2} \, ds\, d\sigma\, dt < \infty\, .
\end{equation}
 Then $v$ is locally Lipschitz.
\end{thm}

\begin{rem}
\begin{enumerate}

\item[i)] Theorem \ref{teo} was proved by Ignat \cite{Ignat_JFA} for $p\in [1,2]$. Moreover, by standard interpolation we have the inclusion $W^{1/p, p}_{div} (\Omega, \Ss^1) \subset W^{1/q, q} (\Omega, \Ss^1)$ for any $p<q$ (the target of such maps being $\Ss^1$, they are always in $L^\infty$).  

\item[(ii)] Note that the $W^{1/p,p}$ assumption naturally excludes ``jump-singularities'' but allows ``oscillations''. For instance, if $p=2$, then the function $\varphi:(-\frac 1 2, \frac 1 2)\to \RR$ defined as $\varphi(x_1)=\log|\log |x_1||$ for $x_1\neq 0$ belongs to $H^{1/2}((-\frac 1 2, \frac 1 2)) = W^{1/2,2} ((-\frac 1 2, \frac 1 2))$. So, setting $\Omega=(-\frac 1 2, \frac 1 2)^2\subset \RR^2$, then the function $u(x_1, x_2):=e^{i\varphi(x_1)}$ belongs to $H^{1/2}(\Omega, \Ss^1)$ and obviously, $L=\{(0, x_2)\, :\, x_2\in (-\frac 1 2, \frac 1 2)\}\subset \Omega$ is an "oscillating" line-singularity of $u$. Theorem \ref{teo} excludes, however, this type of behavior exploiting the additional assumption that $u$ is divergence-free.
\end{enumerate}
\end{rem}

One interesting point is the fundamental role played by a commutator estimate from Constantin, E and Titi, which was used in \cite{CET} to prove that $C^{1/3+\eps}$ solutions of the incompressible Euler equations preserve the kinetic energy. Our proof uses a similar argument to show that the $u$ of Theorem \ref{teo} and the $v$ of Theorem \ref{teo_2} both satisfy some additional balance laws. Such laws hold obviously for smooth solutions but are false in general for distributional solutions. The question of which threshold regularity ensures their validity can be surprisingly subtle. In the case of the incompressible Euler's equations a well-known conjecture of Onsager in the theory of turbulence claims that $1/3$ is the critical H\"older exponent for energy conservation: the ``positive side'' of this conjecture was indeed proved in \cite{CET} (see also \cite{Eyink}), whereas the ``negative side'' is still open, although there have been recently many results in that direction (see for instance \cite{BDS, DS2, DS1, Isett}). 

\section{Entropies and kinetic formulations}

\begin{figure}[htbp]
\center
\includegraphics[scale=0.3,
width=0.3\textwidth]{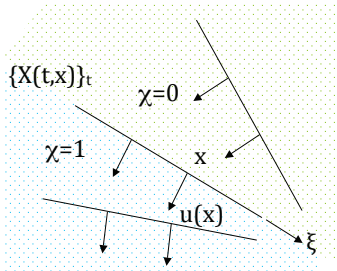} \caption{Characteristics of $u$.} \label{chara}
\end{figure}

\nd The main feature of both problems relies on the concept of characteristic. Assume for the moment that $u$ is a smooth solution of \eqref{contraintes} and fix a point $x\in \Omega$; then the characteristic of $u = \nabla^\perp \psi$ at $x$ is given by 
\be
\label{syst_dyn}
\dot{X}(t,x)=u^\perp(X(t,x))
\ee with the initial condition $X(0,x)=x$. The orbit $\{X(t, x)\}_t$ is a straight line (i.e., $X(t,x)=x+tu^\perp(x)$ for $t$ in some interval around $0$) along which $u$ is perpendicular and constant. A similar conclusion can be drawn for smooth solutions of Burgers' equation, considering the corresponding characteristics, cf. for instance \cite{Dafermos}.
Observe that \eqref{syst_dyn} does not have a direct proper meaning in the case $u\in W^{1/p,p}$ because \`a-priori there is no trace of $u$ defined $\h^1$-a.e. on curves $\{X(t, x)\}_t$. To overcome this difficulty, the following notion of weak characteristic was introduced (see e.g. Jabin-Perthame \cite{Jabin-Perthame}): for every direction $\xi\in \Ss^1$, the function 
$\chi(\cdot, \xi):\Omega\to\{0, 1\}$ is defined as
\be
\label{weak_cara}
\chi(x, \xi)=\begin{cases} 1 &\quad \textrm{ for }\, \, u(x)\cdot \xi>0, \\
0 &\quad \textrm{ for }\, \, u(x)\cdot \xi\leq 0.
\end{cases}\ee
When $u$ is smooth around a point $x\in \Omega$, then for the choice $\xi:=u^\perp(x)$ either $\nabla \chi(\cdot, \xi)$  locally vanishes
(if $u$ is constant in a neighborhood of $x$), or $\nabla \chi(\cdot, \xi)$ is a measure concentrated on the characteristic $\{X(t, x)\}_t$ and oriented by $\xi^\perp$ (see Figure \ref{chara}).
In other words, we have the following ``kinetic formulation" of the problem:
$$\xi\cdot \nabla \chi(x, \xi)=0.$$
Note that the knowledge of $\chi(\cdot, \xi)$ in every direction $\xi\in \Ss^1$ determines completely the vector field $u$ due to the averaging formula 
\be
\label{aver_form}
u(x)=\frac 1 2 \int_{\Ss^1} \xi \chi(x, \xi)\, d\xi  \quad \textrm{ for a.e. }\, \, x\in \Omega.\ee

A similar approach can be used to capture the corresponding characteristics for solutions of Burgers' equation and in fact the work of \cite{Jabin-Perthame} originated from ideas applied first in the theory of scalar conservation laws: inspired by the classical work of Kru\v{z}kov, cf. \cite[Section 6.2]{Dafermos}, a similar ``kinetic formulation'' was introduced first by Lions, Perthame and Tadmor in \cite{LPT} for 
entropy solutions of scalar conservation laws (in any dimension).

\medskip

The key point in the proof of Theorem \ref{teo} consists in showing an appropriate ``kinetic formulation'' for $\wup$-vector fields. 
Indeed Theorem \ref{teo} follows from the following Proposition via an argument of Jabin, Otto and Perthame \cite{JOP02}.

\begin{pro}[Kinetic formulation] 
\label{kinet}
Let $u\in \wup$ with $p\in [1,3]$. For every direction $\xi\in \Ss^1$, the function $\chi(\cdot, \xi)$ defined at \eqref{weak_cara} satisfies the following kinetic equation: 
\be
\label{eqkine}
\xi\cdot \nabla \chi(\cdot, \xi)=0 \quad \textrm{in} \quad {\cal D}'(\Omega). \ee
\end{pro}

\begin{rem}
\label{remh12}
\begin{enumerate}
\item[i)] In Ignat \cite{Ignat_JFA}, the above result was proved for $p\in [1,2]$ and it was conjectured that \eqref{eqkine} still holds for any $p>2$. Proposition \ref{kinet} answers partially to that question for the case $p\leq 3$.

\item[ii)] A ``kinetic averaging lemma'' (see e.g. Golse-Lions-Perthame-Sentis \cite{Golse}) shows that a measurable vector-field $u:\Omega\to \Ss^1$ satisfying \eqref{eqkine} belongs to $H^{1/2}_{loc}$ (due to \eqref{aver_form}). This property can be read as the converse of Proposition \ref{kinet} for the case $u\in H^{1/2}(\Omega, \Ss^1)$. A-posteriori, such a vector field has stronger regularity since it shares the structure described in Theorem \ref{teo}.
\end{enumerate}
\end{rem}

The main concept that is hidden in the kinetic formulation \eqref{eqkine} is that of entropy coming from scalar conservation laws. Indeed, 
for each direction $\xi\in \Ss^1$ we introduce the maps
$\Phi^\xi:\Ss^1\to \RR^2$ defined by
\be
\label{element}
\Phi^\xi(z):=\begin{cases} \xi &\quad \textrm{ for }\, \, z\in \Ss^1,  \, z\cdot \xi>0, \\
0 &\quad \textrm{ for }\, \, z\in \Ss^1, \, z\cdot \xi\leq 0,
\end{cases}\ee
which will be called "elementary entropies". Clearly
$$\Phi^\xi(u(x))=\xi \chi(x, \xi) \quad \textrm{ for a.e. } \quad x\in \Omega$$ 
and \eqref{eqkine} can be regarded as a vanishing entropy production:
$$\nabla \cdot [ \Phi^\xi(u)]=\xi\cdot \nabla \chi(\cdot, \xi)=0 \quad  \textrm{ in } \quad {\cal D}'(\Omega).$$

\medskip

\nd The link between \eqref{contraintes} and scalar conservation laws is the following.
If $u$ is a solution of \eqref{contraintes} of the form $u=(v,h(v))$ (for the flux $h(v)=\pm\sqrt{1-v^2}$) then the divergence-free constraint turns into the 
scalar conservation law 
\begin{equation}
\label{conserv}
v_t + (h(v))_s=0\, .
\end{equation}
From the theory of scalar conservation laws, it is known that, when $h$ is not linear, there is in general no global smooth solution of the Cauchy problem associated to \eqref{conserv}. This leads naturally to consider weak (distributional) solutions of \eqref{conserv} but in this class there are often infinitely many solutions for the same initial data. The concept of entropy solution restores uniqueness, together with good approximation properties with suitable regularizations (see Kru{\v{z}}kov \cite{Kru}). To clarify this notion we recall that an 
entropy - entropy flux pair for \eqref{conserv} is a couple of scalar (Lipschitz) functions $(\eta, q)$ such that $\frac{dq}{dv}=\frac{dh}{dv} \frac{d\eta}{dv}$, which entails that
every smooth solution $v$ of \eqref{conserv} satisfies the balance law $(\eta(v))_t + (q(v))_s=0$.
A solution $v$ of \eqref{conserv} (in the sense of distributions) is called entropy solution if for every convex entropy $\eta$, 
the entropy production $(\eta(v))_t + (q(v))_s$ is a nonpositive measure. We summarize all these
concepts in the following definition for the particular case of Burgers' equation: 

\begin{df}\label{d:entropy_Burg}
An entropy - entropy flux pair $(\eta, q)$ for \eqref{e:Burgers} consists of two (locally) Lipschitz functions $(\eta, q):\mathbb R \to \mathbb R^2$ such that $q' (w) = w \eta' (w)$ for a.e. $w\in \mathbb R$. A distributional solution $v\in L^\infty_{loc}(\Omega)$ of \eqref{e:Burgers} is an
{\em entropy} solution if $(\eta (v))_t + (q(v))_s \leq 0$ for every such pair $(\eta, q)$ with $\eta$ convex. 
\end{df}

The main point of Theorem \ref{teo_2} is to show that $W^{1/p,p}$ {\em weak} solutions of Burgers' equation are in fact
{\em entropy} solutions. \footnote{Heuristically, the link between \eqref{contraintes} and \eqref{e:Burgers} can be understood by approximating $h(v)=-\sqrt{1-v^2}=-1+\frac{v^2}{2}+O(v^4)$ for small $v$ in \eqref{conserv}. Therefore, the link between Theorems \ref{teo} and \ref{teo_2} is the following: in the framework of Theorem \ref{teo_2}, if $\psi$ is a function with
$\psi_t=1-\frac{v^2}{2}$ and $\psi_s=v$, then $\psi$ is a $C^{1,1}$ viscosity solution of the Hamilton-Jacobi equation $\psi_t+\frac{(\psi_s)^2}{2}=1$. Obviously, in the approximation $v$ taken very small, the last equation approximates the eikonal equation $|\nabla \psi|=1$.}

\begin{pro}[Entropy solutions]\label{p:entropy}
Let $v$ be as in Theorem \ref{teo_2}. Then $v$ is a (locally bounded) {\em entropy solution} and moreover
\begin{equation}\label{e:no_shocks}
\left(\frac{v^2}{2}\right)_t + \left(\frac{v^3}{3}\right)_s = 0.
\end{equation} 
\end{pro}

Indeed we will focus in showing only the identity \eqref{e:no_shocks}, since it implies that $v$ is an entropy solution by
\cite[Theorem 2.4]{DOW2} (see also \cite{Panov}). 

In the case of Burgers' (or more generally for conservation laws $v_t + (h(v))_x =0$ with a uniformly convex $h$), entropy solutions $v$ are functions of bounded variation by Oleinik's estimate (see \cite{Dafermos}). The
chain rule of Volpert (cf. \cite[Theorem 3.99]{AFP}) shows then that the entropy production measure $\mu:= (\eta (v))_t + (q (v))_x$  
concentrates on lines (corresponding to "shocks" of $v$): in fact we can use such chain rule to show that \eqref{e:no_shocks} rules out the existence of shocks and then Theorem \ref{teo_2} can be concluded from the classical theory of hyperbolic conservation laws, cf. \cite[Section 11.3]{Dafermos}.
Alternatively we could argue as for Theorem \ref{teo} using the corresponding kinetic formulation, as it is done in \cite[Proposition 3.3]{COW}.

The link between \eqref{contraintes} and \eqref{conserv} suggests to use quantities similar to the entropy - entropy flux pairs $(\eta, q)$ to detect "local" line-singularities of $u$. This idea, which we will explain in a moment, has been used when dealing with reduced models in micromagnetics, e.g., Jin-Kohn \cite{JK00}, Aviles-Giga~\cite{AG99},
DeSimone-Kohn-M\"uller-Otto \cite{DKMO01}, Ambrosio-DeLellis-Mantegazza \cite{AdLM99}, Alouges-Riviere-Serfaty~\cite{ARS02}, Ignat-Merlet \cite{IM09}, \cite{IMpre}, Ignat-Moser \cite{Ignat_Moser}. However in these cases the corresponding entropy production measures usually change sign.
This raised the question of proving the concentration of the entropy production measures on $1$-dimensional sets for those weak solutions with entropy productions which are {\em signed} Radon measures. Partial results are available, see \cite{AKLR,DOW1,DR}, but the general problem is still widely open.

\medskip
 
\nd In the sequel we will always use the following notion of entropy introduced in~\cite{DKMO01} for solutions of the eikonal equation (see also \cite{DLO03, IMpre, JK00}). It corresponds to the entropy - entropy flux pair from the scalar conservation laws, but here the pair is defined in terms of the couple $(v, h(v))$ and not only on $v$.
\begin{df} [DKMO \cite{DKMO01}]
\label{defentrop}
We will say that $\Phi\in C^\infty(\Ss^1,\RR^2)$ is an entropy if
\begin{equation}
\label{condentrop}
\frac{d}{d\theta}\Phi(z)\cdot z \ =\ 0, \quad \textrm{for every $z=e^{i\theta}=(\cos\theta,\sin\theta)\in \Ss^1$.}
\end{equation}
Here, $\frac{d}{d\theta}\Phi(z):=\frac{d}{d\theta}[\Phi(e^{i\theta})]$ stands for the angular derivative of $\Phi$. The set of all entropies is denoted by $\, ENT$.
\end{df}

\medskip

\nd The following two characterizations of entropies are proved in \cite{DKMO01}:

\begin{enumerate}
\item[1.] A map $\Phi\in C^\infty(\Ss^1,\RR^2)$ is an entropy if and only every $u\in C^\infty(\Omega, \RR^2)$ as in \eqref{contraintes} has no entropy production: 
\be  \label{propentrop}
\ncd [\Phi(u)]=0 \quad \textrm{in} \quad {\cal D}'(\Omega). \ee 

\item[2.] A map $\Phi\in C^\infty(\Ss^1,\RR^2)$ is an entropy if and only if there exists a (unique) $2\pi$-periodic function $\f\in C^\infty(\RR)$ such that  for every $z=e^{i\theta}\in \Ss^1$,  
\begin{equation}
\label{prop_ent}
\Phi(z)=\f(\theta)z+ \frac{d\f}{d\theta}(\theta)z^\perp. 
\end{equation}
In this case, 
\be
\label{prop_ent2}
\frac{d}{d\theta}\Phi(z)=\gamma(\theta) z^\perp,\ee where $\gamma \in C^\infty(\RR)$ is the $2\pi$-periodic function defined by 
$ \lambda=\f+\frac{d^2}{d\theta^2}\f$ in $\RR$. 

\end{enumerate}
 
\medskip 

\nd As shown in Ignat-Merlet \cite{IMpre}, these properties can be extended to nonsmooth entropies, in particular to the special class of elementary entropies $\Phi^\xi$ of \eqref{element}, which are maps of bounded variations. Although $\Phi^\xi$ is not a smooth entropy (in fact, $\Phi^\xi$ has a jump at the points $\pm \xi^\perp\in \Ss^1$), the equality \eqref{condentrop} trivially holds in ${\cal D}'(\Ss^1)$.
Moreover, as shown in \cite{DKMO01}, there exists a sequence of smooth entropies $\{\Phi_k\}\subset ENT$ such that
$\{\Phi_k\}$ is uniformly bounded and $\lim_k \Phi_k(z)=\Phi^\xi(z)$ for every $z\in \Ss^1$ (this approximation result follows via \eqref{prop_ent}). 
Therefore, in order to have the kinetic formulation in Proposition \ref{kinet}, we will prove the following result:

\begin{pro} 
\label{pro_equi}
Let $\Phi\in C^\infty(\Ss^1,\RR^2)$ be an entropy. Then for every $u\in \wup$, $p\in [1, 3]$, the identity \eqref{propentrop} holds true.
\end{pro}

\medskip

\nd Note that this result represents an extension to the class of $W^{1/p, p}$-vector fields of the characterization \eqref{propentrop} of an entropy.

\section{Proofs of Proposition \ref{p:entropy} and Proposition \ref{pro_equi}}

Proposition \ref{pro_equi} was proved in \cite{Ignat_JFA} (see also Ignat \cite{Confluentes}) for $p\in [1,2]$ using a duality argument that cannot be adapted to the case $p>2$. We will present the strategy used in \cite{Confluentes} for the case $p=2$, together with a very elementary argument for $p=1$ (cf. Steps 4 and 5 in the proof below) and then we will present a new method
 that enables to conclude in the case $p\leq 3$. However the easier cases $p\in (1,2]$ can be conclude directly from the latter (cf. Step 7 in the proof below).
 
 \medskip
 
\nd {\bf Proof of Proposition \ref{pro_equi}.} Let $\Phi\in C^\infty(\Ss^1,\RR^2)$ be an entropy, i.e., \eqref{condentrop} holds. Let $B\subset\subset \Omega$ be a ball inside and $\{\rho_\eps\}_{\eps>0}$ be a family of standard mollifiers in $\RR^2$ of the form
$$\rho_\eps(x)=\frac 1 {\eps^2} \rho\left(\frac x \eps\right)$$ with $\rho:\RR^2\to \RR_+$ smooth, $\int_{\RR^2} \rho(x)\, dx=1$ and $\supp \rho \subset B_1$ where $B_1$ is the unit ball in $\RR^2$. For $\eps>0$ small enough, we consider the approximation of $u\in \wup$ in $B$ by convolution with $\rho_\eps$:
$$u_\eps=u\star \rho_\eps \quad \textrm{ in } \quad B.$$ Then $u_\eps\in C^\infty(B, \RR^2)$, $\nabla \cdot u_\eps=0$ and $|u_\eps|\leq 1$ in $B$.

\medskip

\nd {\it Step 1. Extension $\tilde \Phi$ of the entropy $\Phi$ to $\RR^2$.} We extend the entropy $\Phi$ to a ``generalized" entropy $\tilde \Phi$ on $\RR^2$. For that, we consider a smooth function $\eta:[0, \infty)\to \RR$ such that $\eta=0$ on $[0, 1/2]\cup [2, \infty)$ and $\eta(1)=1$ and define $\tilde 
\Phi\in C_c^\infty(\RR^2, \RR^2)$ by 
$$
\tilde \Phi(z):=\eta(|z|) \Phi(\frac{z}{|z|}) \, \, \textrm{ for every } \, \, z\in \RR^2\setminus\{0\}.$$ By \eqref{condentrop}, we have that 
\be
\label{ref_extension}
z\cdot D\tilde \Phi(z)z^\perp=|z|z\cdot \frac{\partial \tilde \Phi}{\partial \theta}(z)=
|z|\eta(|z|)z\cdot \frac{d \Phi}{d \theta}(\frac{z}{|z|})\stackrel{\eqref{condentrop}}{=}0, \quad z\in \RR^2,\ee
with the usual notation $(D\tilde \Phi)_{i,j}=\frac{\partial \tilde \Phi_i}{\partial x_j}$.

\medskip

\nd {\it Step 2. Decomposition of $D \tilde \Phi$.} We show that there exist $\Psi \in C^\infty_c( \R^2,\R^2)$
and $\gamma \in C^\infty_c(\R^2,\R)$ such that 
$$
D\tilde \Phi (z)= -2\Psi(z) \otimes z + \gamma(z) Id\quad \textrm{for every } \, \, z\in \RR^2,$$
where $Id$ is the identity matrix (see \cite{DKMO01}). Indeed, one considers 
$$\gamma(z)=\frac{z^\perp\cdot D\tilde \Phi(z)z^\perp}{|z|^2}\quad \textrm{ and }\quad \Psi(z)=\frac{-D\tilde \Phi(z)z+\gamma(z)z}{2|z|^2}, \quad z\in \RR^2.
$$
(Here, $\gamma$ is indeed an extension to the whole plane $\RR^2$ of the function given in \eqref{prop_ent2}.) Denoting $\rv=\frac{z}{|z|}$ and 
$\uv=\frac{z^\perp}{|z|}$ for $z\neq 0$, one checks, using the spectral decomposition, that
\begin{align*}
D\tilde \Phi (z)-\gamma(z)Id&=\bigg(D\tilde \Phi (z)\rv-\gamma(z)\rv\bigg)\otimes \rv + 
\underbrace{\bigg(D\tilde \Phi (z)\uv-\gamma(z)\uv\bigg)}_{=0 \, \, \textrm{by}\, \, \eqref{ref_extension}}\otimes \, \uv=-2\Psi(z) \otimes z \quad \forall z\neq 0.
\end{align*}

\medskip

\nd {\it Step 3. The entropy production $\nabla\cdot [\Phi (u_\eps)]$.} 
For the smooth approximation $u_\eps$, we obtain the entropy production (as in \cite{DKMO01}):
\begin{align}
\nonumber
\nabla\cdot  [\tilde \Phi (u_\eps)]&={\rm Tr} \bigg(D\tilde \Phi(u_\eps)Du_\eps\bigg)\stackrel{Step \, \, 2}{=}-2{\rm Tr} \bigg(\Psi(u_\eps)\otimes u_\eps \, Du_\eps\bigg)+\gamma(u_\eps)
\underbrace{\nabla \cdot u_\eps}_{=0}\\
\nonumber
&=-2\Psi(u_\eps) \cdot  (Du_\eps)^T u_\eps=-\Psi(u_\eps) \cdot \nabla |u_\eps|^2\\
\label{decomp_eps}
&=\Psi(u_\eps) \cdot \nabla \big(1-|u_\eps|^2\big) \quad \textrm{in} \quad B. 
\end{align}

\medskip

\nd {\it Step 4. Proof of \eqref{propentrop} for $p=1$.} The final issue consists in passing to the limit in \eqref{decomp_eps} as $\eps\to 0$.
On one hand, the chain rule implies that $\tilde \Phi(u_\eps)\to \tilde \Phi(u)=\Phi(u)$ in $W^{1,1}(B)$, in particular, 
\be
\label{p=1}
\nabla \cdot [\tilde \Phi(u_\eps)]\to \nabla \cdot [\Phi(u)] \textrm{ in } L^1(B).\ee On the other hand, the chain rule leads to $1-|u_\eps|^2\to 1-|u|^2=0$ in   
$W^{1,1}(B)$, in particular, $$\nabla (1-|u_\eps|^2) \to 0 \textrm{ in } L^1(B).$$ Since $\{\Psi(u_\eps)\}$ is uniformly bounded, the duality $<\cdot, \cdot>_{L^{\infty}(B), L^1(B)}$ leads to  
$$\Psi(u_\eps)\cdot \nabla (1-|u_\eps|^2) \to 0 \textrm{ in } L^1(B),$$ which by \eqref{decomp_eps} and \eqref{p=1} yield $\nabla\cdot [\Phi(u)]=0$ (in $L^1(B)$).

\medskip

\nd {\it Step 5. Proof of \eqref{propentrop} for $p=2$.}  We repeat the above argument using the duality $$<\cdot, \cdot>_{\hb^{-1/2}(B), H^{1/2}_{00}(B)}$$ where $\hb^{-1/2}(B)$ is the dual space of $H^{1/2}_{00}(B)$: 
$$H_{00}^{1/2}(B)=
\{\zeta\in H^{1/2}(B)\,: \, \int_B \int_B \frac{|\zeta(x)-\zeta(y)|^2}{|x-y|^{3}}\, dxdy+\int_B \frac{|\zeta(x)|^2}{d(x)}\, dx<\infty\}$$
with $d(x)={\rm dist}(x, \partial B)$. In fact, $H_{00}^{1/2}(B)$ can be seen as the closure of $C_c^\infty(B)$ in $H^{1/2}(\RR^2)$ (see e.g. \cite{Ignat_JFA} for more details). More precisely, on one hand, the chain rule implies that $\tilde \Phi(u_\eps)\to \tilde \Phi(u)=\Phi(u)$ in $H^{1/2}(B)$, in particular, 
\be
\label{p=2}
\nabla \cdot [\tilde \Phi(u_\eps)]\to \nabla \cdot [\Phi(u)] \textrm{ in } \hb^{-1/2}(B).\ee On the other hand, the chain rule leads to $1-|u_\eps|^2\to 1-|u|^2=0$ in   
$H^{1/2}(B)$, in particular, $$\nabla (1-|u_\eps|^2) \to 0 \textrm{ in } \hb^{-1/2}(B).$$ Since $\Psi(u_\eps)\to \Psi(u)$ in $H^{1/2}(B)$, we conclude that for every $\zeta\in C^\infty_c(B)$, 
$$<\nabla (1-|u_\eps|^2), \zeta \Psi(u_\eps)>_{\hb^{-1/2}(B), H^{1/2}_{00}(B)} \to 0,$$ which by \eqref{decomp_eps} and \eqref{p=2} yield 
$$<\nabla\cdot [\Phi(u)], \zeta>_{\hb^{-1/2}(B), H^{1/2}_{00}(B)}=0.$$ Hence, $\nabla\cdot [\Phi(u)]=0$ in ${\cal D}'(B)$. 

\medskip

\nd {\it Step 6. Proof of \eqref{propentrop} for $p=3$.}  In this case, we use the estimate of Constantin, E and Titi, cf. \cite{CET}.  Let $\zeta\in C^\infty_c(B)$.
By \eqref{decomp_eps}, we write:
\begin{align*}
\int_{B} \zeta(x) \nabla\cdot  [\tilde \Phi (u_\eps)]\, dx&=\int_{B} \zeta(x) \Psi(u_\eps) \cdot \nabla \big(1-|u_\eps|^2\big)\, dx\\
&=\underbrace{\int_{B} \zeta(x) \nabla\cdot \big[\Psi(u_\eps) (1-|u_\eps|^2)\big]\, dx}_{=I_\eps}-\underbrace{\int_{B} \zeta(x) (1-|u_\eps|^2) \nabla\cdot [\Psi(u_\eps)] \, dx}_{=II_\eps}.
\end{align*}

\nd {\it Passing to the limit for $I_\eps$ as $\eps\to 0$.} By dominated convergence theorem, we have that $\Psi(u_\eps) (1-|u_\eps|^2)\to 0$ in $L^1(B)$ so that, after integrating by parts, we conclude $I_\eps\to 0$ as $\eps\to 0$.\\

\medskip

\nd {\it Passing to the limit for $II_\eps$ as $\eps\to 0$.} This part is subdivided in three more steps.

\nd {\it (i)} First, we write for $x\in B$ and  for small $\eps$:
\begin{align*}
1-|u_\eps(x)|^2&=|u|^2\star \rho_\eps(x)-|u\star \rho_\eps(x)|^2\\
&=\int_{\RR^2} |u(x-z)|^2 \rho_\eps(z)\, dz-\bigg(\int_{\RR^2} u(x-z) \rho_\eps(z)\, dz\bigg)\cdot \bigg(\int_{\RR^2} u(x-w) \rho_\eps(w)\, dw\bigg)\\
&=\int_{\RR^2} \int_{\RR^2} u(x-z)\cdot(u(x-z)-u(x-w)) \rho_\eps(z)  \rho_\eps(w)\, dz\, dw\\
&\stackrel{z:=w, \, w:=z}{=}\frac 1 2 \int_{\RR^2}  \int_{\RR^2} \big|u(x-z)-u(x-w)\big|^2 \rho_\eps(z)  \rho_\eps(w)\, dz\, dw\\
&\leq 2\int_{\RR^2} \big|u(x-z)-u(x)\big|^2 \rho_\eps(z) \, dz\\
&\leq \frac{2\|\rho\|_{L^\infty}} {\eps^2} \int_{B_\eps} \big|u(x-z)-u(x)\big|^2 \, dz,
\end{align*}
where we used the inequality $\frac 1 2 |u(x-z)-u(x-w)\big|^2\leq |u(x-z)-u(x)|^2+|u(x-w)-u(x)|^2$ and the properties of the mollifiers, i.e., $\supp \rho_\eps\subset B_\eps$ (that is the ball of radius $\eps$ centered at the origin) and $\int_{B_\eps} \rho_\eps(z)\, dz=1$. 

\nd {\it (ii)} Second, we write the last term in $II_\eps$ as $\nabla\cdot [\Psi(u_\eps)]={\rm Tr} \bigg(D\Psi(u_\eps)\nabla u_\eps\bigg)$. Moreover, since 
$\int_{B_\eps} \partial_j \rho(\frac z \eps)\, dz=0$ for $j=1,2$, we observe that 
\begin{align*}
\partial_j u_\eps(x)&=u\star \partial_j  \rho_\eps(x)=\frac 1 {\eps^3} \int_{B_\eps}u(x-z) \partial_j  \rho(\frac{z} \eps)\, dz=\frac 1 {\eps^3} \int_{B_\eps}\big(u(x-z)-u(x)) \partial_j  
\rho(\frac{z} \eps)\, dz\\
&\leq \frac{\|\nabla \rho\|_{L^\infty}} {\eps^3} \int_{B_\eps} \big|u(x-z)-u(x)\big| \, dz,
\end{align*}
for $j=1,2$.

\nd {\it (iii)} Third, using Jensen's inequality, we deduce by {\it (i)} and {\it (ii)}:
\begin{align}
|II_\eps|&\leq  \frac C \eps \int_B \Big(\int_{B_\eps} \hspace{-6mm}-\, \, \, |u(x-z)-u(x)|\, dz\Big) \Big(\int_{B_\eps} \hspace{-6mm}-\, \, \, |u(x-z)-u(x)|^2\, dz\Big) \, dx \nonumber\\
&\leq \frac C \eps \int_B  \Big(\int_{B_\eps} \hspace{-6mm}-\, \, \, |u(x-z)-u(x)|^3\, dz\Big)^{1/3} 
 \Big(\int_{B_\eps} \hspace{-6mm}-\, \, \, |u(x-z)-u(x)|^3\, dz\Big)^{2/3}\, dx \nonumber\\
 &= \frac C \eps \int_B \int_{B_\eps} \hspace{-6mm}-\, \, \, |u(x-z)-u(x)|^3\, dz\, dx \\
 &= \frac C {\eps^3} \int_B \int_{B_\eps} |u(x-z)-u(x)|^3\, dz\, dx \nonumber\\
&\stackrel{|z|\leq \eps}{\leq} \int_B \int_{B_\eps} \frac{|u(x-z)-u(x)|^3}{|z|^3}\, dz\, dx = \int_B \int_{B_\eps (x)} \frac{|u(x)-u(y)|^3}{|y-x|^3}\, dy\, dx\label{e:vanishing}\, .
\end{align}
Since $u\in W^{1/3, 3}(B)$, the integral
\[
\int_{B\times B} \frac{|u(x)-u(y)|^3}{|y-x|^3}\, dy\, dx
\]
is finite and thus the last integral in \eqref{e:vanishing} converges to $0$ as $\eps\downarrow 0$.
Therefore, we conclude that \eqref{propentrop} holds for $p=3$.
\medskip

\nd {\it Step 7. Proof of \eqref{propentrop} for $p\in (1,3)$.} By
Gagliardo-Nirenberg embedding: $L^\infty\cap W^{1/p, p}\subset W^{1/3,3}$ (see \cite{BBM_lifting}, Lemma D.1) and thus, one concludes by Step 6.

\medskip

\nd Since $B\subset \Omega$ was an arbitrarily chosen ball, \eqref{propentrop} follows in $\Omega$. \qed

\begin{proof}{ of Proposition \ref{p:entropy}}
We use computations very similar to those of Step 6 in the previous proof to show that \eqref{e:no_shocks} holds. More precisely
consider a family of standard mollifiers $\rho_\eps$, but this time in the space variable $s$ only: $\rho \in C^\infty_c (]-1,1[)$ and $\rho_\eps (s) = \frac{1}{\eps} \rho (\frac{s}{\eps})$. We still use the notation $v_\eps = v\star \rho_\eps$ for the convolution of $v$ and $\rho$ in the space variable {\em only}, namely
\[
v\star \rho_\eps (t,s) = \int v (t, s-\sigma) \rho_\eps (\sigma)\, d\sigma\, .
\]
Fix a smooth test function $\zeta\in C^\infty_c (\Omega)$. Our goal is to show that
\begin{align}\label{e:regularized}
\lim_{\eps\downarrow 0} \underbrace{\int_\Omega \left(\textstyle{\frac{v_\eps^2}{2}} \zeta_t + \textstyle{\frac{v_\eps^3}{3}} \zeta_x\right)}_{=:J_\eps} = 0\, .
\end{align}
This in turn would imply that \eqref{e:no_shocks} holds and the Proposition would then follow from \cite[Theorem 2.4]{DOW2}. Observe that, although
we are only mollifying in space, we can conclude from \eqref{e:Burgers} that
\begin{equation}\label{e:Burgers_mollified}
(v_\eps)_t + \left(\textstyle{\frac{v^2 \star \rho_\eps}{2}}\right)_s = 0 \qquad \mbox{in $\Omega_\eps = \{(s,t)\in \Omega : {\rm dist}\, ((s,t), \partial \Omega) > \eps\}$.}
\end{equation}
In particular, for $\eps$ sufficiently small, $v_\eps$ turns out to be $C^1$ on the support of $\zeta$.
Integrating by parts, using the chain rule and then subtracting \eqref{e:Burgers} we easily reach
\begin{align}
J_\eps &= - \int v_\eps \zeta \left((v_\eps)_t + \left(\textstyle{\frac{v_\eps^2}{2}}\right)_s\right) \stackrel{\eqref{e:Burgers_mollified}}{=} - \frac{1}{2} \int v_\eps \zeta (v_\eps^2 - v^2\star\rho_\eps)_s\nonumber\\
&= \frac{1}{2} \underbrace{\int (v_\eps)_s \zeta (v_\eps^2 - v^2 \star \rho_\eps)}_{=: I_\eps} + \frac{1}{2} \int v_\eps \zeta_s (v_\eps^2 - v^2\star \rho_\eps)\, .\label{e:identita}
\end{align}
Observe that the second integral in \eqref{e:identita} goes to $0$ because $v_\eps$ is uniformly bounded in $L^3$ (indeed by assumption it is bounded in $L^4$) and $v_\eps^2 - v^2\star \rho_\eps$ converges to $0$ strongly in $L^{3/2}$ (in fact by assumption it converges even in $L^2$). We thus need to show that $I_\eps$ converges to $0$. Following the same computations of the Steps 6 and 7 in the previous proof we can easily show that:
\begin{align*}
|(v_\eps)_s (t,s)| &= \frac{2}{\eps} \left|\int_{-\eps}^\eps \hspace{-6mm}-\, \, \,  (v (t,s -\sigma) - v(t,s)) \rho' \left(\frac{\sigma}{\eps}\right)\, d\sigma\right|\nonumber\\ 
&\leq \frac{C}{\eps} \left(\int_{-\eps}^\eps \hspace{-6mm}-\, \, \,  |v (t, s-\sigma) - v(t,s)|^3\, d\sigma\right)^{1/3}\\
|v_\eps^2 - v^2*\rho_\eps| (t,s) &= \frac{1}{2\eps^2} \left|\int\int (v (t, s-\sigma) - v (t,s-\sigma'))^2 \rho \left(\frac{\sigma}{\eps}\right)
\rho \left(\frac{\sigma'}{\eps}\right)\, d\sigma d\sigma'\right|\nonumber\\
&\leq {C} \int_{-\eps}^\eps \hspace{-6mm}-\, \, \, |(v(t, s-\sigma)- v(t,s)|^2\, d\sigma\\
&\leq C \left(\int_{-\eps}^\eps \hspace{-6mm}-\, \, \,  |v (t, s-\sigma) - v(t,s)|^3\, d\sigma\right)^{2/3}\, .
\end{align*}
Recalling that ${\rm supp}\, (\zeta) \subset I \times K \subset\subset I \times J$ for some closed interval $K$, we conclude
\begin{align*}
|I_\eps| &\leq \frac{C}{\eps^2} \int_I \int_K \int_{-\eps}^\eps |v(t, s-\sigma) - v (t,s)|^3\, d\sigma\, ds\, dt\\
&\leq C\int_I \int_K \int_{s-\eps}^{s+\eps} \frac{|v(t,s)-v(t,\sigma)|^3}{|s-\sigma|^2}\, d\sigma\, ds\, dt\, .
\end{align*}
Since by assumption
\[
\int_I \int_{K\times K}  \frac{|v(t,s)-v(t,\sigma)|^3}{|s-\sigma|^2}\, d\sigma\, ds\, dt < \infty\, ,
\]
we obviously conclude that $I_\eps\to 0$.
\qed
\end{proof}

\section{Proofs of Proposition \ref{kinet}, Theorem \ref{teo} and Theorem \ref{teo_2}}

\begin{proof}{ of Proposition \ref{kinet}} For every $\xi\in \Ss^1$, 
the non-smooth "elementary entropies" $\Phi^\xi:\Ss^1\to \RR^2$ given by \eqref{element}
can be approximated by a sequence of smooth entropies $\{\Phi_k\}\subset ENT$ such that
$\{\Phi_k\}$ is uniformly bounded and with $\lim_k \Phi_k(z)=\Phi^\xi(z)$ for every $z\in \Ss^1$. Indeed, this smoothing result follows by \eqref{prop_ent}: if one writes $\xi=e^{i\theta_0}$ with $\theta_0\in(-\pi, \pi]$, then the unique 
$2\pi$-periodic function $\f\in C(\RR)$ satisfying \eqref{prop_ent} for $\Phi^\xi$ is given by:
$$\f(\theta)=\xi\cdot z {\bf 1}_{\{z\cdot\xi>0\}}=\cos(\theta-\theta_0) {\bf 1}_{\{\theta-\theta_0\in (-\pi/2, \pi/2)\}} \quad \textrm{ for }\, \, z=e^{i\theta}, \, \theta\in (-\pi+\theta_0, \pi+\theta_0).$$ 
By \eqref{prop_ent} for $\Phi^\xi$, the choice of $\f'$ is fixed at the jump points $\pm \xi^\perp\in \Ss^1$: 
$$\f'(\theta)=-\sin (\theta-\theta_0) {\bf 1}_{\{\theta-\theta_0\in (-\pi/2, \pi/2)\}} \quad \textrm{ for}\, \,  \theta\in (-\pi+\theta_0, \pi+\theta_0).$$
Now, one regularizes $\f$ by $2\pi-$periodic functions $\f_k\in C^\infty(\RR)$ that are uniformly bounded in $W^{1, \infty}(\RR)$ and $\lim_k \f_k(\theta)=\f(\theta)$ as well as
 $\lim_k \f'_k(\theta)=\f'(\theta)$ for every $\theta\in \RR$. Thus, the desired (smooth) approximating entropies $\Phi_k$ are given by $\f_k$ via \eqref{prop_ent}.
Therefore, Proposition \ref{pro_equi} implies that for every {$u\in \wup$ (with $p\in [1,3]$)},
one has $\int_\Omega \Phi_k(u)\cdot \nabla \zeta\, dx=0$ for every $\zeta\in C^\infty_c(\Omega)$ and by the dominated convergence theorem, we pass to the limit $k\to \infty$ and conclude that $$0=\ncd [\Phi^\xi(u)]=\ncd [\xi \chi(\cdot, \xi)]=\xi\cdot \nabla \chi(\cdot, \xi) \quad \textrm{in} \quad {\cal D}'(\Omega).$$
\qed
\end{proof}

\proof{ of Theorem \ref{teo}}
It is a consequence of Proposition \ref{kinet} combined with the strategy of Jabin-Otto-Perthame (see Theorem 1.3 in \cite{JOP02}). For completeness of the writing, let us recall the main steps of that argument:
let $u:\Omega\to \Ss^1$ be a measurable function that satisfies \eqref{eqkine} for every $\xi\in \Ss^1$. Notice that the divergence-free condition is automatically satisfied (in ${\cal D}'(\Omega)$) because
of \eqref{aver_form}. The first step consists in defining a $L^\infty$-trace of $u$ on each segment $\Sigma\subset \Omega$. More precisely, if $\Sigma:=\{0\}\times [-1,1]\subset \Omega$, then there exists a trace
$\tilde u\in L^\infty(\Sigma, \Ss^1)$ such that
$$\lim_{r\to 0} \frac 1 r \int_{-r}^r \int_{-1}^1 |u(x_1, x_2)-\tilde u(x_2)|\, dx_2 dx_1=0$$ and for each Lebesgue point $(0,x_2)\in \Sigma$ of $u$, one has $u(0,x_2)=\tilde u(x_2)$. Observe that this step is straightforward in the case of $u\in \wun$; however, it is essential for example in the case of $p>1$. The second step is to prove that if the trace $\tilde u$ of $u$ is orthogonal at $\Sigma$ at some point, then $\tilde u$ is almost everywhere orthogonal at $\Sigma$ (which coincides with the classical principle of characteristics for smooth vector fields $u$). The key point for that resides in a relation of order of characteristics of $u$, i.e., for every two Lebesgue points $x, y\in \Omega$ of $u$ with the segment $[x,y]\subset \Omega$, the following implication holds: 
$$u(x)\cdot (y-x)>0 \Rightarrow u(y)\cdot (y-x)>0.$$ The final step consists is proving that on any open convex subset $\omega \subset \Omega$ with $d=\dist(\omega, \partial \Omega)>0$, only two situations may occur: either two 
characteristics of $u$
intersect at $P\in \Omega$ with $\dist(P, \omega)<d$ and $u(x)=\pm \frac{(x-P)^\perp}{|x-P|}$ for $x\in \omega\setminus \{P\}$, or $u$ is $1/d$-Lipschitz in $\omega$, i.e.,
$$|u(x)-u(y)|\leq \frac 1 d |x-y|, \quad \textrm{ for every } \, x,y\in \omega.$$
(In this last case, every two characteristics passing through $\omega$ may intersect only at distances $\geq d$ outside $\omega$). Note that $u$ may have infinitely many vortex points $P_k$ and any vortex point has 
degree one, but the orientation $\alpha_k$ of the vortex point $P_k$ could change or not in $\Omega$.  
\qed

\begin{proof}{ of Theorem \ref{teo_2}} As shown in Proposition \ref{p:entropy}, $v$ is an entropy solution. As such, we conclude from the classical Oleinik's estimate (cf. \cite[Theorem 11.2.1]{Dafermos}) that $v_x$ is a Radon measure and hence that $v$ is in fact $L^\infty_{loc}$ and $BV_{loc}$. On the other hand the equality \eqref{e:no_shocks} implies that $v$ is shock-free in $\Omega$ (cf. for instance the proof of \cite[Corollary 2.5]{DOW2}). In particular it follows from \cite[Theorem 11.3.2]{Dafermos} that $v$ is everywhere continuous and therefore from \cite[Theorem 11.3.5]{Dafermos} that it is locally Lipschitz.
\qed
\end{proof}

\bibliographystyle{plain}

\end{document}